\documentclass{amsart}
\usepackage{amssymb}
\usepackage{subfig}
\usepackage{epsfig}
\usepackage{graphicx}
\usepackage{epic}
\usepackage{float}
\usepackage[mathscr]{eucal}
\usepackage{xspace}
\usepackage[latin1]{inputenc}
\usepackage{amssymb,amsgen,amsbsy,amsopn,amsfonts,graphicx}
\usepackage[dvips]{color}
\usepackage{multicol}
\usepackage{mathrsfs}

\newcommand{\R}{\mathbb{R}}

\newcommand{\Z}{\mathbb{Z}}
\newcommand{\N}{\mathbb{N}}

\newcommand{\LL}{\mathcal{L}}

\newcommand{\po}{\partial}

\newcommand{\ve}{\varepsilon}

\newcommand{\loc}{{\text{loc}}}
\newcommand{\X}{\times}

\renewcommand{\d}{\delta}
\renewcommand{\l}{\lambda}

\renewcommand{\b}{\beta}

\renewcommand{\k}{\kappa}
\newcommand{\z}{\zeta}

\newcommand{\om}{\omega}

\newcommand{\M}{{\mathcal M}}
\newcommand{\V}{{\mathcal V}}

\newcommand{\rbf}{\mathbf{r}}

\renewcommand{\subset}{\subseteq}

\newcommand{\TV}{\text{TV}\,}

\newcommand{\Rf}{{\mathfrak R}}

\newcommand{\per}{\operatorname{per}}

\theoremstyle{plain}
\newtheorem{theorem}{Theorem}[section]

\theoremstyle{definition}

\theoremstyle{remark}

\numberwithin{equation}{section}

\usepackage{amsthm}

\begin{document}

\title[Bakhvalov systems with damping]{ An extension of Bakhvalov's theorem for systems of conservation laws with damping}

\author{Hermano Frid}
\address{Instituto de Matem\'atica Pura e Aplicada - IMPA\\ Estrada Dona Castorina, 110\\
Rio de Janeiro, RJ, 22460-320, Brazil}
\email{hermano@impa.br}
\thanks{}

 
\begin{abstract} For $2\X2$ systems of conservation laws satisfying Bakhvalov conditions, we present  a class of damping terms that still yield  the existence of global solutions with periodic initial data of possibly large bounded total variation per period. We also address the question of the decay of the periodic solution. As applications we consider the systems of isentropic gas dynamics, with pressure obeying a $\gamma$-law, for the physical range $\gamma\ge1$, and also for the ``non-physical'' range $0<\gamma<1$,  both in the classical Lagrangian and Eulerian formulation, and in the relativistic setting. We give complete details for the case $\gamma=1$, and also analyze the general case when $|\gamma-1|$ is small. Further, our main result also establishes the decay of the periodic solution.        

\end{abstract}

\maketitle

\section{Introduction}\label{S:1}

In \cite{Ni}, Nishida shows that for the $2\X 2$ system of compressible isotermic gas dynamics, where the pressure satisfies a $\gamma$-law with $\gamma=1$,
it is possible to construct a global solution using the Glimm method \cite{Gl} for any initial data of bounded  total variation, taking values in the physical region, with no smallness restriction on the value of the total variation.   Following Nishida's work, Bakhvalov, in the important paper \cite{Ba}, establishes general conditions for a $2\X2$ system of conservation laws to enjoy the property of allowing the construction of a global solution by the Glimm method with no restriction of smallness in the value of the initial total variation.     
The purpose of this paper is to establish a general result extending the theorem of Bakhvalov, in \cite{Ba}, for $2\X2$ systems of conservation laws systems with  damping. More specifically, we prescribe a general form for the damping term, with 4 degrees of freedom, and,  besides Bakhvalov's conditions in \cite{Ba}, we impose additional conditions concerning our prescribed family of damping terms so as to obtain general conditions for the existence of a  global solution, whose initial data may have large total variation, depending only on how big is the region where the original Bakhvalov's conditions are satisfied.  Our main result also includes the decay of the periodic entropy solution of the system with damping.  We then first apply the general result to the $2\X2$ system of isotermic gas dynamics, that is, pressure given by a so called $\gamma$-law with $\gamma=1$,  in Lagrangian and Eulerian coordinates, as well as the relativistic version. Finally, we also discuss the application of the main result to the general isentropic gas dynamics with pressure obeying a $\gamma$-law, when $\gamma$ is sufficiently close to 1.  In what follows we assume that the reader  has some  basic knowledge of the theory of conservation  laws as presented, e.g.,  in \cite{Sm, Se, Ho, Da}.

So, let us consider a $2\X2$ system of conservation laws in the form 
\begin{equation}\label{e1.1}
\begin{cases}
u_t+f_1(u,v)_x+g_1(u,v)=0,\\
v_t+f_2(u,v)_x+g_2(u,v)=0,
\end{cases}
\end{equation}
or in the form
 \begin{equation}\label{e1.1'}
\begin{cases}
U_1(u,v)_t+f_1(u,v)_x+ \tilde g_1(u,v)=0,\\
U_2(u,v)_t+f_2(u,v)_x+\tilde g_2(u,v)=0,
\end{cases}
\end{equation}
where $(u,v)\mapsto (U_1,U_2)$ is  locally a change of coordinates,
and the initial data 
\begin{equation}\label{e1.2}
(u(x,0),v(x,0))=(u_0(x),v_0(x)).
\end{equation}

Let $\om(U)$ and $\z(U)$ be a pair of Riemann invariants for \eqref{e1.1} (or \eqref{e1.1'}), $U=(u,v)$, defined for $U\in\V$, with $\V$ a domain in $\R^2$. Let $\bar U$ be given by $(\om,\z)(\bar U)=(0,0)$.

We recall that \eqref{e1.1} (or \eqref{e1.1'}) is said to satisfy Bakhvalov's conditions (see \cite{Ba})  on a domain $\V$ if the following are satisfied, where $\z=L_i(\om;\om_r,\z_r)$ denotes  $i$-left shock curve of the points that can be connected on the left by an $i$-shock to $(\om_r,\z_r)$, and  $\z=R_i(\om;\om_r,\z_r)$ denotes  $i$-right shock curve of the points that can be connected on the right by an $i$-shock to $(\om_r,\z_r)$, $i=1,2$:
\begin{eqnarray}
B_1: && \sup_{i,\V}{|\lambda_i(\om,\z)|}<\infty. \nonumber \\
B_2: && \forall (\om,\z)\in\V,\; 0<\frac{\partial R_1}{\partial \om},\,\frac{\partial L_1}{\partial \om}<1,\;1<\frac{\partial R_2}{\partial \om},\,\frac{\partial L_2}{\partial \om}<+\infty,\;w\neq w_0
\nonumber \\
B_3: && \mbox{ If $\z_r=R_i(\om_r;\om_l,\z_l),\,i=1,2,$ then shock curves $\z = R_i(\om;\om_l,\z_l),\,\om\le \om_l$}
 \nonumber \\
    && \mbox{ and $\z = L_i(\om;\om_r,\z_r),\,\om \ge \om_r$ intersect only in points $(\om_l,\z_l),\,(\om_r,\z_r)$.} \nonumber \\
B_4: && \mbox{If four points $(\om_l,\z_l),\,(\om_r,\z_r),\,(\om_m,\z_m)$ and $(\hat{\om}_m,\hat{\z}_m)$ satisfy}
\nonumber \\
    && \mbox{ $\z_m = R_2(\om_m;\om_l,\z_l),\, \z_r = R_1(\om_r;\om_m,\z_m)$, $\hat{\z}_m = R_1(\hat{\om}_m;\om_l,\z_l)$ and } \nonumber \\
     && \z_r = R_2(\om_r;\hat{\om}_m,\hat{\z}_m), \text{then } (\om_l - \hat{\om}_m) + (\hat{\z}_m - \z_r)\leq(\z_l - \z_m)+(\om_m - \om_r).\nonumber
\end{eqnarray}

 We make the following  assumptions:
\begin{enumerate}
\item[{\bf(D1)}] We assume that \eqref{e1.1} (or \eqref{e1.1'}) is a strictly hyperbolic genuinely nonlinear system and that, with respect to the pair $(\om,\z)$, \eqref{e1.1} (or \eqref{e1.1'}) satisfies the Bakhvalov conditions $B_1$--$B_4$  on $\V$.  These, in particular, ensure  that both $\om$ and $\z$ decrease across the shocks of both families on $\V$, and , if $\Rf(x,t;U_L,U_M)$ and $\Rf(x,t;U_M,U_R)$ are the solutions of the Riemann problems formed for the two pairs of constant states $(U_L,U_M)$ and $(U_M,U_R)$, respectively,  and $\Rf(x,t;U_L,U_R)$ is the solution of the Riemann problem for the pair $(U_L,U_R)$, then
\begin{equation}\label{e1.3}
\LL(\Rf(U_L,U_R))\le \LL(\Rf(U_L,U_M))+  \LL(\Rf(U_M,U_R)),
\end{equation}
where 
\begin{equation}\label{e1.4}
\LL(\Rf(x,t;U_L,U_R)):= \sum_{1-\text{shocks}} |[\om(\Rf)]| +\sum_{2-\text{shocks}} |[\z(\Rf)]|,
\end{equation}
and $|[\om(\Rf)]|=\om(U_{\text{left}})-\om(U_{\text{right}})$, $U_{\text{left}}$ being the state to the left of the shock wave, and $U_{\text{right}}$ the one to the right.

\item[{\bf(D2)}] In the case of system  \eqref{e1.1}, we assume that $G(U)=(g_1(U),g_2(U))$, satisfies 
\begin{equation}\label{e1.5}
G(U)= a\om\rbf_1+  b\z \rbf_2,
\end{equation}
where $a>0, b>0$  and $\rbf_1,\rbf_2$ are normalized right-eigenvectors such that 
\begin{equation}\label{e1.5'}
\rbf_1\cdot\nabla \om=1,\quad \rbf_1\cdot\nabla \z=0,\quad \rbf_2\cdot\nabla\om=0,\quad \rbf_2\cdot\nabla\z=1.
\end{equation}
In the case of system \eqref{e1.1'}, denoting
\begin{equation}\label{e1.1''}
A(u,v)=\left[\begin{matrix} \frac{\po U_1}{\po u} & \frac{\po U_1}{\po v}\\ \frac{\po U_2}{\po u} &\frac{\po U_2}{\po v} \end{matrix}\right],
\end{equation}
we then have that $\tilde G(u,v)=(\tilde g_1(u,v),\tilde g_2(u,v))$ is defined by the formula
\begin{equation}\label{e1.1'''}
\tilde G(u,v):=A(u,v) G(u,v),
\end{equation}
where  $G(u,v)$ is given by \eqref{e1.5} and now $\rbf_1, \rbf_2, \om, \z$ are eigenvectors and Riemann invariants for the matrix $A^{-1}(U)\nabla f(U)$, with $f=(f_1,f_2)$.

\item[{\bf(D3)}] There exists $\d_0>0$, such that, for $0<\d<\d_0$, if $U_L,U_R$ are connected by a 1-shock wave (or a 2-whock wave),  $\Lambda:=\left[\begin{matrix} a&0\\0&b\end{matrix}\right]$, with $a>0,b>0$, as  in the previous item,  the Riemann solution (in the plane of Riemann invariants) with left state
$\tilde Z_L:=e^{-\Lambda \d}Z_L$ and right state $\tilde Z_R:=e^{-\Lambda \d}Z_R$, with $Z_L:=(\om_L,\z_L)$, $Z_R:=(\om_R,\z_R)$,  $\om_L=\om(U_L)$, etc., satisfies
\begin{multline}\label{e1.4'}
\LL(\Rf(\tilde Z_L,\tilde Z_R))\le \LL(\Rf(Z_L,Z_R))=\om(U_L)-\om(U_R)\\ \text{($=\z(U_L)-\z(U_R)$, in the case of a  2-shock)}.
\end{multline}

\item[{\bf(D4)}] The system \eqref{e1.1} admits a strictly convex entropy $\eta_*$ satisfying $\eta_*(\bar U)=0$, $\nabla\eta_*(\bar U)=0$, and
\begin{equation}\label{e1.5''}
\nabla_{{}_U}\eta_*(U)\cdot G(U)=a\,\om\po_\om\eta_{*}(U(\om,\z))+b\,\z\po_\z\eta_{*}(U(\om,\z))\ge0,
\end{equation}
 which, when $a=b$, simply means that  $\eta_*$ is radially nondecreasing in the Riemann invariants $(\om,\z)$-plane. We recall that this means that there exists a companion function $q_*$, the entropy flux, such that $\nabla\eta_*\nabla f=\nabla q_*$.

\end{enumerate}

We remark that, from {\bf(D2)}, we see, in particular, that our family of damping terms, given by \eqref{e1.5} (or \eqref{e1.1'''}) has 4 degrees of freedom: 2 corresponding to the choice of $\bar U$, and 2 corresponding to the choice of $a$ and $b$ in \eqref{e1.5}.

We assume that $(u_0,v_0)$ is a periodic function, with period, say, 1, with local bounded total variation, and we set $\b_0:=\TV_{\text{per}}(u_0,v_0)$, where $\TV_{\per}$ stands for total variation in $x$ over one period.  Let 
$$
\mu_0:=\int_0^1\eta_*(U_0(x))\,dx,
$$
$$
\Sigma_0=\{U\in\R^2\,:\, \eta_*(U)\le \mu_0\},
$$ 
and
$$
\Sigma_{0,R}:=\{U\in\R^2\,:\,\operatorname{dist}(U;\Sigma_0)<R\}.
$$

Let $\V_0\subset\R^2$ be the invariant domain
$$
\V_0:=\{U\in\R^2\,:\,\om(U)\le \sup_{x\in[0,1]}\om(U_0(x)),\ \z(U)\ge \inf_{x\in[0,1]}\z(U_0(x))\}.
$$

\begin{theorem}\label{T:1.1}
There exists a constant $C_0>0$, depending only on the flux functions in \eqref{e1.1} (or \eqref{e1.1'}), such that if  
$$
\Sigma_{0,C_0 \b_0}\cap \V_0\subset \V,
$$ 
the problem \eqref{e1.1},\eqref{e1.2} (or \eqref{e1.1'},\eqref{e1.2}) possesses a global entropy solution $U(x,t)$, with locally bounded total variation, such that 
\begin{equation}\label{e1.T1}
U(x,t)\in \Sigma_{0,C_0\b_0}\cap \V_0, \quad\text{for all $(x,t)\in\R\X[0,\infty)$}
\end{equation}
Moreover, if $\eta_*$ is strictly dissipative in the sense that, except at $U=\bar U$, the strict inequality holds in \eqref{e1.5''}, we have 
\begin{equation}\label{e1.T1'}
\lim_{t\to\infty}\int_0^1|U(x,t)-\bar U|\,dx=0.
\end{equation}

\end{theorem}  

The proof of Theorem~\ref{T:1.1} is given in Section~\ref{S:2}.  In Section~\ref{S:3}, we apply Theorem~\ref{T:1.1} to the $2\X 2$ system of compressible esoteric gas dynamics, in which the pressure is given by a $\gamma$-law with $\gamma=1$. Finally,  in Section~\ref{S:4}, we address the application of Theorem~\ref{T:1.1} to the case of a  general $\gamma$-law pressure with $\gamma$ sufficiently close to 1.

\section{Proof of Theorem~\ref{T:1.1}}\label{S:2}

In this section we prove Theorem~\ref{T:1.1}. 

We construct an approximate solution for \eqref{e1.1},\eqref{e1.2} by using a fractional step procedure where, for each time interval $[nh,(n+1)h)$, we make two successive iterations in the following way. Assume the approximate solution has been defined at $t=nh$, such that $U^h(x,nh)=U^h((m+a_n)l,nh)$, for $x\in ((m-\frac12)l,(m+\frac12)l)$,  for some randomly chosen number $a_n\in(-\frac12,\frac12)$, where $m\in\Z$, and  $l=1/N_l$, for some $N_l\in\N$, satisfies the usual CFL-condition 
$$
\frac{l}{h}\ge  \sup_{U\in\V}\{|\l_i(U)|\,:\, i=1,2\}.
$$
 As a first step, we define  $\hat U^h(x,t)$, for $(x,t)\in\R\X[nh,(n+1)h)$, according to  Glimm's method, that is, defining it over the mesh rectangle $((m-\frac12)l,(m+\frac12)l)\X[nh,(n+1)h)$ as the Riemann solution centered at $(ml,nh)$ with left state $ U^h((m+a_n)l,nh)$ and right state $U^h((m+1+a_n)l,nh)$.  In the second step, we define $U^h(x,t)$, for $nh\le t<(n+1)h$, $x\in\R$, as the solution of
\begin{equation}\label{e1.9}
 \begin{cases}
\frac{dU}{dt}(x,t)=-G(U(x,t)), \quad nh<t<(n+1)h\\
U(x,nh)=\hat U^h((m+a_{n+1})l,(n+1)h-0), \quad  \text{for $(m-\frac12)l<x<(m+\frac12)l$},
 \end{cases}
 \end{equation}
again for some randomly chosen number $a_{n+1}\in (-\frac12,\frac12)$.  We close the recursive definition of $U^h(x,t)$ by setting 
 $$
 U^h(x,(n+1)h)=U^h((m+a_{n+1})l,(n+1)h-0),\quad \text{for $(m-\frac12)l<x<(m+\frac12)l$}.
 $$
Let us also define
$$
\hat U(x,(n+1)h)=\hat U^h((m+a_{n+1})l,(n+1)h-0),\quad \text{for $(m-\frac12)l<x<(m+\frac12)l$},
$$
which would be the prescription for constructing the Glimm solution, for the corresponding homogeneous system.

{}From the above construction, we see that, in each time interval $t\in[nh,(n+1)h)$, $U^h(x,t)$ is piecewise constant in $x$, for each fixed $t$.   We define the functional
 \begin{equation}\label{e1.10}
\LL(U^h(x,nh)):= \sum_{1-\text{shocks}}^{\per} [\om(U^h)] +\sum_{2-\text{shocks}}^{\per} [\z(U^h)],
\end{equation}
where $ \sum\limits_{1-\text{shocks}}^{\per}$ means that the sum is over all 1-shocks over one period resulting from the solution of the Riemann problems for each of the discontinuities over one period and  $[\om(U)]=\om(U_{\text{left}})-\om(U_{\text{right}})$, with a similar definition for  $\sum\limits_{2-\text{shocks}}^{\per}$ and  $[\z(U)]$. 

We recall that Bakhvalov's theorem implies that
\begin{equation}\label{e1.11}
\LL(\hat U^h(x,(n+1)h))\le \LL(U^h(x,nh)).
\end{equation}
To get the desired inequality
\begin{equation}\label{e1.12}
\LL(U^h(x,(n+1)h)\le \LL(U^h(x,nh)),
\end{equation}
it then suffices to prove
\begin{equation}\label{e1.13}
\LL(U^h(x,(n+1)h)\le \LL(\hat U^h(x,(n+1)h)).
\end{equation}
Observe that, from \eqref{e1.9}, in passing from $\hat U^h$ to $U^h$, the Riemann invariants satisfy 
\begin{equation}\label{e1.14}
\begin{aligned} 
&\frac{d\om}{dt}=-a\om,\\
&\frac{d\z}{dt}=-b\z.
\end{aligned}
\end{equation}
To obtain inequality \eqref{e1.13}, we first observe that condition {\bf(D3)} implies that inequality \eqref{e1.4'} also holds if $U_L$ and $U_4$ are any two constant states, not necessarily connected by a $1$-shock or a $2$-shock, as a consequence of Bakhvalov condition, where we use the fact that if $P$ and $Q$ are connected by a rarefaction wave of the first or the second family, then $e^{-\Lambda\d}P$ and $e^{-\Lambda\d}Q$ are also connected by a rarefaction wave of the same family. Therefore, inequality \eqref{e1.13} follows from the validity of the corresponding inequality for each individual discontinuity and its corresponding transformation by the application of $e^{-\Lambda h}$ to both points forming the discontinuity, written in the Riemann invariants coordinates.  This suffices to prove the uniform boundedness of the total variation per period for $U^h(x,t)$. 
 
 We now briefly explain how the control of the $L^\infty$ norm of $U^h$ is achieved. For that we need to use the strictly convex entropy $\eta_*$ whose existence is assumed. 
Because of condition {\bf(D4)}, we have 
\begin{align}\label{e1.15}
\eta_*(\int_0^1 U^h(x,t)\,dx)&\le \int_0^1\eta_*(U^h(x,t))\,dx\le \int_0^1\eta_*(U_0^h(x))\,dx\\ 
&+\sum_{j=1}^{j=[t/h]} \int_0^1(\eta_*(U^h(x,jh-0))-\eta_*(U^h(x,jh+0)))\,dx.\notag
\end{align}
Observe that, as proven in \cite{Lx}, the summation at the end of \eqref{e1.15} is nonpositive in the limit as $h\to0$ and it is always bounded by $C\b_0 t$, for some constant $C>0$. Since $\eta_*$ is strictly convex, this means that the mean value
$\int_0^1 U^h(x,t)\,dx$ is always inside the region $\Sigma_0$, in a sufficiently small time interval. Since, as we have just proven,  the total variation per period is bounded at any time by $\frac12 C_0\b_0$, for some constant $C_0$ only depending on the system, we conclude that $U^h (x,t)$ assumes its values in $\Sigma_{0,C_0\b_0}$. This argument should be performed in a stepwise manner in intervals of time of a fixed length $T$, extracting subsequences for an appropriate choice of the random sequence $\{a_n\}$   by means of a diagonal argument, as in \cite{Fr,FP}.  

Now, the proof of the consistence of the above  Glimm scheme with fractional step follows by standards arguments whose central point relies on the original consistence argument in \cite{Gl}, which then implies the convergence of $U^h$ to a bounded function with bounded total variation per period $U$ which is a weak solution of \eqref{e1.1'} (or \eqref{e1.1'}, as the case may be). Also, by the same argument in \cite{Lx} , we obtain that $U$ satisfies the entropy inequality
\begin{equation}\label{e1.16}
\int_{\R\X(0,\infty)}\{ \eta_*(U)\varphi_t+q_*(U)\varphi_x\}\,dx\,dt+\int_{\R}\eta_*(U_0(x))\varphi(x,0)\,dx\ge0,
\end{equation}
for any $0\le\varphi\in C_0^\infty(\R^2)$.  

To prove the decay property \eqref{e1.T1} we follow the approach in \cite{CF} based on the compactness of the scaling sequence $U^\ve(x,t)=U(\ve^{-1}x,\ve^{-1}t)$ as $\ve\to0$. We notice that $U^\ve$ is an entropy weak solution of a system like \eqref{e1.1} (or \eqref{e1.1'}) only that instead of $G(U)$ (resp., $\tilde G(U)$), we now have $\ve^{-1} G(U^\ve)$ (resp., $\ve^{-1}\tilde G(U^\ve)$).  Clearly, inequality \eqref{e1.16} is also satisfied by $U^\ve$ and this implies, through a standard argument in the theory of distributions,  that 
$$
\eta_*(U^\ve)_t+q_*(U^\ve)_x \in \{\text{bounded subset of $\M_\loc^-(\R^2\X(0,\infty))$}\},
$$
where  $\M_\loc^-(\R^2\X(0,\infty))$ denotes the space of nonpositive Radon measures of locally bounded total variation. Now, under the assumption that $(\eta_*,q_*)$ is strictly dissipative, this implies that for any entropy-entropy flux pair $(\eta,q)$ for the system \eqref{e1.1} (or \eqref{e1.1'}), satisfying $\eta(\bar U)=0$, $\nabla\eta (\bar U)=0$, we have 
$$
\eta(U^\ve)_t+q(U^\ve)_x \in \{\text{bounded subset of $\M_\loc(\R^2\X(0,\infty))$}\},
$$  
where  $\M_\loc(\R^2\X(0,\infty))$ denotes the space of Radon measures of locally bounded total variation.  Indeed, this follows from the fact that, under the strictly dissipative condition,  given any such pair $(\eta,q)$ we may find a constant $C_\eta$ such that $\eta+C_\eta\eta_*$ satisfies the dissipative condition \eqref{e1.5''} on a given compact domain. Here we also use the strict convexity of $\eta_*$ to handle the dissipative inequality around $\bar U$. Therefore, through a well know interpolation result (see, {\em e.g.}, \cite{Mu, Ta}) it follows that 
$$
\eta(U^\ve)_t+q(U^\ve)_x \in \{\text{compact subset of $H_\loc^{-1}(\R^2\X(0,\infty))$}\}.
$$  
This allows us to apply the DiPerna's compactness theorem for general $2\X2$ strictly hyperbolic genuinely nonlinear systems in \cite{Di}. Hence, we deduce that $U^\ve$ is compact in $L_\loc^1(\R\X(0,1))$. Now, taking any convergent subsequence, still denoted by $U^\ve$, we see that, if $U^\ve\to\tilde U$ in $L_\loc^1(\R\X(0,\infty))$, we must have
$$
G(\tilde U(x,t))=0,\quad \text{(or $\tilde G(\tilde U(x,t))$, in the case of \eqref{e1.1'})},
$$
for a.e.\ $(x,t)\in\R\X(0,\infty)$, which implies that $\tilde U(x,t)\equiv \bar U$, a.e. In particular, we deduce that the whole sequence $U^\ve$ converges to $\bar U$, a.e. Now using the same argument as in \cite{CF} we finally conclude that \eqref{e1.T1'} holds.

\section{Application to isotermic gas dynamics} \label{S:3}

In this section we present, as an application of Theorem~\ref{T:1.1}, the system of isothermic gas dynamics, first in classical mechanics, both in Lagrangian and Eulerian coordinates, and also in relativistic mechanics.   
 
We consider first the case where the corresponding homogeneous system (i.e., without damping)  is the $p$-system of isothermic gas dynamics
\begin{equation}\label{e1.6'}
\begin{cases}
u_t-v_x=0,\\
v_t+p(u)_x=0.
\end{cases}
\end{equation}
Here, $u$ represents the specific volume,  $v$ represents the velocity, and the pressure is given by $p(u)=\k u^{-\gamma}$,  with $\k>0$, which is the so called $\gamma$-law for the pressure.  Here we will consider the case of an isotermic gas which corresponds to the adiabatic exponent $\gamma=1$. For simplicity we take $\k=1$. 
We take $\V=\{(u,v)\in\R^2\,:\, u>0\}$. In this case, for $\bar U=(\bar u,\bar v)$, with $\bar u>0$,  the Riemann invariant are the classical ones
\begin{equation}\label{e1.7}
\om=v-\bar v+\log \frac{u}{\bar u},\quad \z=v-\bar v-\log\frac{u}{\bar u},
\end{equation}
Therefore, in this case, we have
\begin{equation}\label{e1.8}
\rbf_1= \frac{1}{2}(u,1),\quad \rbf_2= \frac{1}{2}(-u,1).
\end{equation}
Hence, if $a=b$,  $G$ has the form
$$
\begin{aligned}
G(U)&= \frac{a}{2} ((v-\bar v+\log\frac{u}{\bar u})(u,1)+(v-\bar v-\log\frac{u}{\bar u})(-u,1))\\
         &=a(u\log\frac{u}{\bar u}, v-\bar v).
\end{aligned}         
$$
Finally, we easily verify that the well known convex entropy
$$
\eta(u,v)= \frac12 (v-\bar v)^2-\log \frac{u}{\bar u}+\frac1{\bar u}(u-\bar u),
$$
with entropy flux $q(u,v)=(v-\bar v)(1/u- 1/\bar u)$, satisfies condition {\bf (D4)} and it is clearly strictly dissipative in the case $a=b$.  Actually, when $a=b$, the strict dissipation follows from the strict convexity since in this case \eqref{e1.5''} amounts to require that $\eta$ is strictly increasing with respect to the radius in a polar coordinates system centered at $\bar U$.   Thus, for $\gamma=1$ and $a=b$, \eqref{e1.6'} with damping becomes
\begin{equation}\label{e3.1}
\begin{cases}
u_t-v_x=a\,u\log\frac{u}{\bar u},\\
v_t+(\frac1{u})_x=a(v-\bar v).
\end{cases}
\end{equation}

System \eqref{e1.6'} is also frequently presented in the so called Eulerian coordinates, where it  takes a form like \eqref{e1.1'}.
 In this case, we get the following system
\begin{equation}\label{e1.6''}
\begin{cases}
\rho_t+(\rho v)_x=0,\\
(\rho v)_t+(\rho v^2+p(\rho))_x=0.
\end{cases}
\end{equation}
Here, $\rho$ represents the density, $v$ is again the velocity, and the pressure is given by $p(\rho)=\k \rho^\gamma$. Again we consider the esoteric case where 
$\gamma=1$, and we take $\k=1$.  
We take $\V=\{(\rho,v)\in\R^2\,:\, \rho>0\}$.  Now, with $\bar U=(\bar \rho, \bar v)$,  the expressions for the  classical Riemann invariants take the forms 
\begin{equation}\label{e1.7-2}
\om=v-\bar v+\log\frac{\rho}{\bar \rho},\quad \z=v-\bar v-\log \frac{\rho}{\bar \rho},
\end{equation}
and the right-eigenvectors are given by
\begin{equation}\label{e1.8-2}
\rbf_1= \frac{1}{2}(\rho,1),\quad \rbf_2= \frac{1}{2}(-\rho,1).
\end{equation}
Therefore, if $a=b$, 
$$
G(\rho,v)= a(\rho\log\frac{\rho}{\bar \rho}, v-\bar v).
$$ 
Now, we have $U_1(\rho,v)=\rho$, $U_2(\rho,v)=\rho v$, and so
$$
A(\rho,v)=\left[\begin{matrix} 1 &0\\ v & \rho  \end{matrix}\right].
$$
In this case,  a simple calculation shows that 
$$
\tilde G(\rho,v)=a(\rho\log\frac{\rho}{\bar \rho}, \rho( v-\bar v)+ \rho v\log\frac{\rho}{\bar \rho}).
$$
Also, concerning condition {\bf(D4)},  the system \eqref{e1.6''} admits the following well known convex entropy 
$$
\eta(\rho,v)=\frac12 \rho( v-{\bar v})^2+ \rho\log\frac{\rho}{\bar\rho}-\rho+\bar\rho,
$$
which is clearly strictly dissipative in the case $a=b$, for the already mentioned reason.  Thus, for $a=b$, \eqref{e1.6''} with damping becomes
\begin{equation}\label{e3.2}
\begin{cases}
\rho_t+(\rho v)_x=a\,\rho\log\frac{\rho}{\bar \rho},\\
(\rho v)_t+(\rho v^2+p(\rho))_x=a(\rho( v-\bar v)+ \rho v\log\frac{\rho}{\bar \rho}).
\end{cases}
\end{equation}

Another example of the importance of  having the prescription for the damping for systems in the more general form \eqref{e1.1'} is provided by the relativistic version of the Euler equations of isentropic gas dynamics, namely,
\begin{equation}\label{e1.rel1}
\begin{cases}
 \frac{\po}{\po t} U_1(\rho,v)+\frac{\po}{\po x} U_2(\rho,v)=0,\\
 \frac{\po}{\po t} U_2(\rho,v)+\frac{\po}{\po x} (U_2(\rho,v) v+p(\rho))=0,
 \end{cases}
 \end{equation}
 where 
 \begin{equation}\label{e1.rel2}
 \begin{aligned}
 & U_1=\frac{v}{c^2} U_2+\rho\\
 &U_2=(p+\rho c^2)\frac{v}{c^2-v^2},\\
 &p(\rho)= \k^2 \rho^{\gamma},
 \end{aligned}
 \end{equation}
 where $\gamma\ge1$. Here we just consider the case $\gamma=1$; $\k$ is a positive constant, and $c$ is the speed of light. Also, for simplicity we only write the formulas for $\bar U=(1,0)$. Here we take $\V=\{ (\rho,v)\in\R^2\,:\, \rho>0,\ |v|\le c\}$.  Existence of a global solution in $BV$ of \eqref{e1.rel1}, with $\gamma=1$,  for initial data in $BV$, with no smallness restriction on the total variation, was obtained in \cite{ST}, where it was observed that the shock curves for this system enjoy the same features as those for \eqref{e1.6'} and \eqref{e1.6''}, which allowed Nishida's global existence result for large initial total variation.  
 
 In this case, we have
 $$
 A(\rho,v)=\left[ \begin{matrix} (\k^2+c^2)\frac{v^2}{c^2(c^2-v^2)}+1 & (\k^2+c^2)\rho \frac{2v}{(c^2-v^2)^2}\\  \\ (\k^2+c^2)\frac{v}{c^2-v^2} & (\k^2+c^2)\rho\frac{c^2+v^2}{(c^2-v^2)^2}\end{matrix}\right].
$$
On the other hand, the classical Riemann invariants are given by the expressions
\begin{equation}\label{e1.rel3}
\begin{aligned}
&\om:= \frac{c}{2}\log\frac{c+v}{c-v}+\frac{c^2}{\k^2+c^2}\log\rho,\\
& \z:=\frac{c}{2}\log\frac{c+v}{c-v}-\frac{c^2}{\k^2+c^2}\log\rho.
\end{aligned}
\end{equation}
The Riemann invariants for any $\bar U$ are $\om(U)-\om(\bar U)$, $\z(U)-\z(\bar U)$. For the right-eigenvalues we then have the expressions
\begin{equation}\label{e1.rl4}
\begin{aligned}
&\rbf_1:=\frac12\left(\rho(1+\frac{\k^2}{c^2}),1-\frac{v^2}{c^2}\right)\\
&\rbf_2:=\frac12\left(-\rho(1+\frac{\k^2}{c^2}),1-\frac{v^2}{c^2}\right).
\end{aligned}
\end{equation}
Using the formula
$$
\tilde G(\rho,v):= A(\rho,v) G(\rho,v),
$$
with
$$
G(\rho,v):= a (\zeta \rbf_1+ \om \rbf_2),
$$
we obtain the prescribed formula for the damping in the relativistic isentropic Euler equations of gas dynamics 
\begin{equation}\label{e1.rel5}
\begin{cases}
 \frac{\po}{\po t} U_1(\rho,v)+\frac{\po}{\po x} U_2(\rho,v)+\tilde g_1(\rho,v)=0,\\
 \frac{\po}{\po t} U_2(\rho,v)+\frac{\po}{\po x} (U_2(\rho,v) v+p(\rho)) + \tilde g_2(\rho,v)=0,
 \end{cases}
 \end{equation}
 To obtain a convex entropy for the system \eqref{e1.rel1} we may just apply the general result of Lax in \cite{Lx}, which establishes the existence of such a convex entropy defined on any compact over which  the system is strictly hyperbolic and genuinely nonlinear. The strict convexity implies the strictly dissipative property in condition {\bf(D4)}.

Condition {\bf(D1)}, that is, the fact that the system satisfies Bakhvalov's conditions,  in all the above formulations for the isothermic gas dynamics system, is well known. The key fact to be used here is that the shock curves for these systems, in the Riemann invariants plane, belonging to each of the two characteristic families, are all translations of the same curve based on a fixed state in the Riemann invariants plane. Condition {\bf(D2)}, in all these examples, is satisfied since  the damping terms satisfy the prescribed formula in this condition. 

As for condition {\bf(D3)}, using the property enjoyed by the shock curves just mentioned,  the inequality \eqref{e1.3} can be proved with the help of Fig.1, in the following manner. In Fig.~1 we see a 1-shock from $P_1$ to $P_2$, and the Riemann solution for the discontinuity between $P_1' =e^{-ah}P_1$ and $P_2'=e^{-ah}P_2$.  These are the images of $P_1$ and $P_2$ when we pass from $\hat U^h$ to $U^h$, but to get a better view of what is going on and be able to compare lengths we translate the entire composite wave curve representing the Riemann solution so that $P_1'$ assumes the position of $P_1$. Due to the fact that 1-shock curves are translations of one another, after the mentioned translation, have the coinciding two 1-shock curves as shown in Fig.~1.   By the second Bakhvalov condition, the slope of the shock curves is greater than 0 and less than 1.  We then deduce that
\begin{align*}
\om(P_1)-\om(P_2)&= \om(P_1)-\om(P_3')+\om(P_3')-\om(P_2)\\
&\ge \om(P_1)-\om(P_3')+\z(P_3')-\z(P_2)\\
&\ge \om(P_1)-\om(P_3')+\z(P_3')-\z(P_2'),
\end{align*}
where for the inequality in the second line we use the fact that for a 1-shock curve $\Delta\om\ge\Delta\z$, by Bakhvalov condition $B_2$, and the inequality in the third line is obvious. We also observe, concerning the inequality in the second line, that we are  taking advantage from the fact, peculiar for the $\gamma$-law case $\gamma=1$,  that the shock curves are translations from one another. Similarly, we prove the corresponding  inequality for the 2-shock waves (see Fig.2). In this case, after translating the composite wave curve representing the Riemann solution, to improve our way to compare lengths, we further translate the 2-shock curve in the translated Riemann solution, so as to make it  lie on the original 2-shock curve connecting $P_1$ and $P_2$, and $P_2'$ goes over to $P_3''$, and from $P_3''$ to $P_2'$ we draw the translated 1-shock, which originally connects  $P_1$ and the middle state $P_3'$.  In particular, we have $\z(P_3')-\z(P_2')=\z(P_1)-\z(P_3'')$ and 
$\om(P_1)-\om(P_3')=\om(P_3'')-\om(P_2')$, and so, using the fact that in a 2-shock we have $\Delta\z\ge\Delta \om$,  we get
\begin{align*}
\z(P_1)-\z(P_2)&= \z(P_1)-\z(P_3'')+\z(P_3'')-\z(P_2)\\
&\ge \z(P_1)-\z(P_3'')+\om(P_3'')-\om(P_2)\\
&\ge \z(P_1)-\z(P_3'')+\om(P_3'')-\om(P_2')\\
&=\om(P_1)-\om(P_3')+\z(P_3')-\z(P_2').
\end{align*}

      These facts imply that inequality \eqref{e1.3} holds, as was to be proved.

\section{Application to the general $\gamma$-law for $\gamma$ close to 1}\label{S:4}

In this section we outline the application of Theorem~\ref{T:1.1} to the general $\gamma$-law compressible isentropic gas dynamics ofr $\gamma$ sufficiently close to 1. 

First we recall that in \cite{Ba} it was proven that systems \eqref{e1.6'} and \eqref{e1.6''} with $0<\gamma<1$ satisfy Bakhvalov's conditions $B_1$--$B_4$,  for the usual Riemann invariants, namely, 
\begin{equation}
\begin{aligned}\label{e4.1}
&\om=v-\bar v- \frac{2\gamma^{1/2}\k}{\gamma-1}(u^{-(\gamma-1)/2}-{\bar u}^{-(\gamma-1)/2}) ,\\
&\z=v-\bar v+ \frac{2\gamma^{1/2}\k}{\gamma-1}(u^{-(\gamma-1)/2}-{\bar u}^{-(\gamma-1)/2}),
\end{aligned}
\end{equation}
for \eqref{e1.6'}, while  the corresponding Riemann invariants for \eqref{e1.6''} are obtained by just making $u=1/\rho$, $\bar u=1/\bar \rho$ in \eqref{e4.1}.   

As observed in \cite{FP}, the same is true for the relativistic version  \eqref{e1.rel1}, for $p(\rho)=\k^2\rho^\gamma$. Namely,  when $0<\gamma<1$, Bakhvalov's
conditions  $B_1$--$B_4$ are satisfied for the standard Riemann invariants given by ({\em cf.} \cite{ST}) 
\begin{equation}
\begin{aligned}\label{e4.2}
&\om=\frac12\log\frac{c+v}{c-v}-\frac12\log\frac{c+\bar v}{c-\bar v}+c\int_{\bar\rho}^\rho\frac{\sqrt{p'(s)}}{p(s)+sc^2}\,ds \\
&\z=\frac12\log\frac{c+v}{c-v}-\frac12\log\frac{c+\bar v}{c-\bar v}-c\int_{\bar\rho}^\rho\frac{\sqrt{p'(s)}}{p(s)+sc^2}\,ds.
\end{aligned}
\end{equation}

On the other hand, for $\gamma>1$, DiPerna \cite{DP} introduced a new family of Riemann invariants for \eqref{e1.6'} and \eqref{e1.6''} for which these systems satisfy   Bakhvalov's conditions  $B_1$--$B_4$ over certain region contained in the physical domain $\{\,\rho>0\,\}$ in the Riemann invariants plane, which includes any compact subset of the physical domain, as long as  $\gamma>1$ is sufficiently close to 1. DiPerna's Riemann invariants, when properly normalized, converge locally uniformly with all their derivatives to the standard Riemann invariants for $\gamma=1$ ({\em cf.} \eqref{e1.7}, \eqref{e1.7-2}) as $\gamma\to1+$.  

In \cite{FP} this fact was extended to the relativistic version \eqref{e1.rel1}. Namely, in \cite{FP} it is proven that DiPerna's formula defining the family of new Riemann invariants in terms of the classical Riemann invariants provide also in the relativistic case a pair of Riemann invariants for which Bakhvalov's conditions  $B_1$--$B_4$ are satisfied in a region $\V$ in the physical domain $\{\,\rho>0,\; v^2<c^2\, \}$ which includes any compact in the physical domain, as long as $\gamma>1$ is sufficiently close to 1.

So, concerning condition  {\bf(D1)}, for $0<\gamma<1$, for systems \eqref{e1.6'} or \eqref{e1.6''},  we can take the whole Riemann invariants plane as $\V$, while for the relativistic version \eqref{e1.rel1} we may take as $\V$ the subset of the Riemann invariants plane  corresponding to  physical domain $\{\, v^2<c^2\,\}$. On the other hand, for $1<\gamma<2$, for systems \eqref{e1.6'} or \eqref{e1.6''}, we may take as $\V$ the region of the Riemann invariants plane for which DiPerna's Riemann invariants satisfy Bakhvalov's conditions $B_1$--$B_4$, which includes any compact contained in the physical domain $\{\,\rho>0,\; v^2<c^2\, \}$, as long as $\gamma>1$ is sufficiently close to 1. 

Condition {\bf(D2)} only prescribe the formula for the damping terms. As for {\bf(D4)}, it is satisfied in the general case, $0<\gamma$, when $a=b$. 

Therefore, it only remains to discuss condition {\bf(D3)} for $\b(\gamma):=|\gamma-1|>0$. We keep assuming $a=b$. We claim that condition {\bf(D2)} is also satisfied in any  compact region in the physical domain, for $\b(\gamma)>0$ sufficiently small. So, let $K$ be a given compact region in the physical domain in the Riemann invariants plane. Let also $V\supset K$ be a neighborhood of $K$ with compact closure contained in the physical domain and $\po V\in C^\infty$. We introduce the functions $\Theta_1(\gamma, \d, W_1,W_2)$ and $\Theta_2(\gamma, \d, Z_1,Z_2)$, with $W_1,Z_1\in K$, $W_2,Z_2\in \bar V$, $W_2$ is connected to the right to $W_1$ by a 1-shock curve relatively to the 
$\gamma$-law system, in either classical or relativistic version, $Z_2$ is connected to the right to $Z_1$ by a 2-shock curve also pertaining to the $\gamma$-law system in either  case, and $0<\d\le 1$. $\Theta_1$ is defined as follows. Given $\gamma>0$, $0<\d\le 1$, and $W_1\in K$, for each $W_2$ in the right 1-shock curve pertaining to the $\gamma$-law, say, in the classical version, we consider the Riemann problem between $W_1'=e^{-a\d}W_1$, as a left state, and $W_2'=e^{-a\d}W_2$, as a right state. It is known that this Riemann problem is solved by drawing a 1-shock curve issuing from $W_1'$ and an inverse 2-wave curve issuing from $W_2'$, consisting of the states that can be connected to the left to $W_2'$ by either a 2-shock or a 2-rarefaction wave. Then, we translate the composite wave curve representing the  Riemann solution together with the line connecting $W_1'$ to $W_2'$, which is parallel to that connecting $W_1$ and $W_2$, until the the line connecting $W_1'$ to $W_2'$ lies entirely over the line connecting $W_1$ to $W_2$, with $W_1'$ coinciding with $W_1$. Denote by $W_3'$ the translated middle state in the Riemann solution connecting $W_1'$ to $W_2'$. $\Theta_1(\gamma, \d, W_1,W_2)$ is then defined as the acute angle between the line connecting $W_2$ to $W_3'$ and the vertical axis passing through $W_2$ ({\em cf.} Fig.~\ref{theta}).  As for $\Theta_2(\gamma,\d, Z_1, Z_2)$, similarly, we first solve the Riemann problem connecting $Z_1'=e^{-a\d}Z_1$ and $Z_2'=e^{-a\d}Z_2$, and translate the composite wave curve together with the line connecting $Z_1'$ to $Z_2'$, until it lies over the line connecting $Z_1$ to $Z_2$.  But now, instead of considering the translated Riemann solution and its middle state, we consider the composite formed by the translation of the 2-shock curve ending in $Z_2'$, so that the translated curve starts now at $Z_1$, and the translation of the 1-wave curve, say, 1-shock curve, so that the corresponding translated curve ends up at $Z_2'$ ({\em cf.} Fig.~\ref{fractional2}). Then we consider the middle state of the composite curve obtained after these translations and call it $Z_3''$. $\Theta_2$ is defined as the angle formed between the line connecting $Z_2$ to $Z_3''$ and the horizontal axis passing through $Z_2$.
We complete the definition of both $\Theta_1$ and $\Theta_2$ for $\d=0$, by defining $\Theta_1(\gamma,0,W_1,W_2)$ as the angle that the tangent to 1-shock curve connecting $W_1$ to $W_2$, at $W_2$, for with the vertical axis through $W_2$. Similarly, $\Theta_2(\gamma,0, Z_1,Z_2)$ is the angle that the tangent to the 2-shock curve connecting $Z_1$ to $Z_2$, at $Z_2$, forms with the horizontal axis through $Z_2$.      

For $\gamma=1$, we have $\Theta_1(1,\d, W_1,W_2)>\pi/4$, for any $W_1\in K$, $W_2\in\bar V$ connected by a 1-shock to $W_1$, and $0\le \d\le1$. Assuming the continuity  of $\Theta_1$ with respect $(\gamma,\d, W_1,W_2)$, we deduce that $\Theta_1(\gamma,\d,W_1,W_2)$ for $|\gamma-1|$ sufficiently small, depending on $K$ and $\bar V$. Similarly, assuming the continuity of $\Theta_2$ with respect to $(\gamma,\d,Z_1,Z_2)$,  we obtain that $\Theta_2(\gamma,\d, Z_1,Z_2)>\pi/4$, for $|\gamma-1|$ sufficiently small. Now, as we have seen in the case $\gamma=1$, the fact that $\Theta_1>\pi/4$ and $\Theta_2>\pi/4$ implies the condition {\bf(D3)}, and this finishes the verification of conditions {\bf(D1)}--{\bf(D4)} for $|\gamma-1|>0$ sufficiently small. It would remain only to prove the continuity of $\Theta_1$ and 
$\Theta_2$ with respect to their arguments, which, although intuitive, requires some topological digression and we leave the details to be given elsewhere.

\subsection{Extension to $|a-b|>0$.}\label{SS:4.1}  We now outline how one can extend the verification of {\bf(D1)}--{\bf(D4)} for $|a-b|>0$, sufficiently small.  We keep considering the systems \eqref{e1.6'}, \eqref{e1.6''} and \eqref{e1.rel1}, for $|\gamma-1|>0$. We first see that the conditions {\bf(D1)} and {\bf(D2)} do not require any checking concerning whether  $a=b$ or $a\ne b$. So, only conditions {\bf(D3)} and {\bf(D4)} need to be discussed concerning the extension of the application of Theorem~\ref{T:1.1} to the gas dynamics systems with a damping term according to the prescription in {\bf(D2)} when $|a-b|>0$. Let us set $\mu:=b-a$. We again assume $K$ to be a given compact in the physical domain in the Riemann invariants plane and $V$ to be an open neighborhood of $K$ whose closure is a compact in the physical domain. 

In order to extend {\bf(D3)} to the case where $|a-b|>0$, we consider again the functions $\Theta_1$ and $\Theta_2$ introduced above, but now we consider them also depending on $\mu$, thus $\Theta_1=\Theta_1(\gamma, \d,\mu, W_1,W_2)$, $\Theta_2=\Theta_2(\gamma,\d,\mu,Z_1,Z_2)$, so that  $\Theta_1(\gamma,\d,W_1,W_2)$ and 
$\Theta_2(\gamma,\d,Z_1,Z_2)$ as defined before now correspond to the values of these functions when $\mu=0$. The definition of these functions for $\mu\ne0$ is similar as that for $\mu=0$, but now, in the case of $\Theta_1$, the line connecting $W_1'=e^{-a\d}W_1$ and $W_2'=e^{-b\d}W_2$ is no longer parallel to that connecting $Z_1$ and $Z_2$, which cause no significant difference. We again translate the line  connecting $W_1'$ to $W_2'$, together with the composite wave curve representing the Riemann solution between $W_1'$ and $W_2'$, until $W_1'$ goes over unto $W_1$. Again $\Theta_1$ is defined as the acute angle that the line connecting $W_2$ to the translated middle state of the Riemann solution $W_3'$ makes with the vertical axis passing through $W_2$. Similarly, we extend the definition of $\Theta_2$ from the case $\mu=0$ to the case where $\mu\ne 0$.   

Again we assume continuity of $\Theta_1$ and $\Theta_2$ with respect to $(\gamma,\d,\mu,W_1,W_2)$ and $(\gamma,\d,\mu,Z_1,Z_2)$, respectively. Now,  we have $\Theta_1(1,\d,0,W_1,W_2)>\pi/4$, for all 
$\d\in[0,1]$, $W_1\in K$ and $W_2\in\bar V$ belonging to the right 1-shock curve issuing from $W_1$.  Also, similarly,  $\Theta_2(1,\d,0,Z_1,Z_2)>\pi/4$,  for all 
$\d\in[0,1]$, $Z_1\in K$ and $Z_2\in\bar V$ belonging to the right 2-shock curve issuing from $Z_1$. Therefore, continuity allows us to obtain $\ve_0>0$ sufficiently small such that, for $|\gamma-1|<\ve_0$ and $|\mu|<\ve_0$, we still have   $\Theta_1(\gamma,\d,\mu,W_1,W_2)>\pi/4$ and  $\Theta_2(\gamma,\d,\mu,Z_1,Z_2)>\pi/4$, and this implies the verification of {\bf(D3)}. 

As for {\bf(D4)},  let us consider the strictly convex entropies corresponding to each of the systems \eqref{e1.6'}, \eqref{e1.6''} and \eqref{e1.rel1}  which are strictly dissipative when $a=b$. We write     
\begin{align*}
\nabla_{{}_U}\eta_*(U)\cdot G(U)&=a\left(\om\po_\om\eta_{*}(U(\om,\z))+\z\po_\z\eta_{*}(U(\om,\z))\right) + (b-a)\,\z\po_\z\eta_{*}(U(\om,\z))\\
                                                          &=A+B.
\end{align*}
Now, by the strict convexity of $\eta_*$,  and $\eta_*(\bar U)=0$, $\nabla\eta_*(\bar U)=0$,  we see that for $|b-a|<\ve_0$, with $\ve_0$ as in the preceding paragraph, there exists a sufficiently small neighborhood of $\bar U$, say, for $|U-\bar U|<r$, for $r>0$ sufficiently small, such that $|B|\le A$ for $|U-\bar U|<r$, and so $A+B\ge0$ in this neighborhood. On the other hand, for $|U-\bar U|\ge r$, with $U$ belonging to the compact corresponding to $K$, $\tilde K=(\om,\z)^{-1}(K)$, since $A>0$ on $\tilde K\cap \{|U-\bar U|\ge r\}$, we see that for $|b-a|$ sufficiently small $A>|B|$ on $\tilde K\cap \{|U-\bar U|\ge r\}$. Therefore, we may find $0<\ve_1\le \ve_0$ such that for
$|b-a|<\ve_1$ condition {\bf(D4)} is satisfied.

\section*{Acknowledgement}

The author gratefully acknowledges the support from CNPq, through grant proc.~303950/2009-9, and FAPERJ, through grant E-26/103.019/2011.
This paper is an outcome of the IMPA-Brown University agreement which supported some visits of the author to Brown University where he had the opportunity of enjoying  many discussions about this research with Costas Dafermos to whom he is deeply indebted for many important suggestions and mathematical tips.

 \begin{figure}[t]
  \caption{A 1-shock before and after the second part of the fractional step.} 
 \centering
  \includegraphics{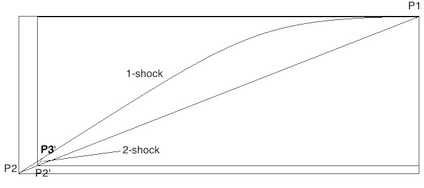} 
  \label{fractional1}
\end{figure}

 \begin{figure}[b]
  \caption{A 2-shock before and after the second part of the fractional step.} 
 \centering
  \includegraphics{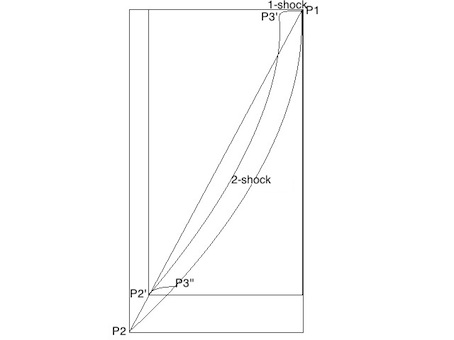} 
  \label{fractional2}
\end{figure}

 \begin{figure}[t]
  \caption{The definition of the function $\Theta_1$.} 
 \centering
  \includegraphics{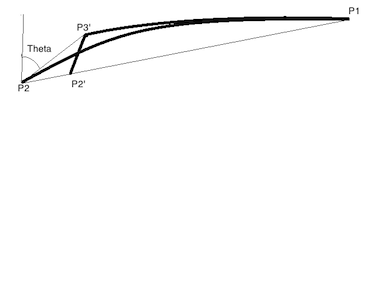} 
  \label{theta}
\end{figure}

\end{document}